\documentclass[oneside,10pt]{amsart}
\usepackage{amsmath,amsfonts, latexsym,amssymb}
\usepackage[active]{srcltx}

\theoremstyle{definition}

\numberwithin{equation}{section}


\begin{document}

\title[Closed  maximal ideals some  Fr\'{e}chet 
algebras \ldots  ]{Closed maximal ideals ideals in  
some Fr\'{e}chet algebras  of holomorphic functions}

\author{Romeo Me\v strovi\' c}

\address{University of Montenegro, Maritime Faculty Kotor,
\,\, Dobrota 36, 85330 Kotor, Montenegro, 
e-mail: {\tt romeo@ucg.ac.me}}

\begin{abstract}
The space $F^p$ ($1<p<\infty$)  consists of  
 all   holomorphic  functions $f$  on the open unit disk $\Bbb D$ such that    
   $ \lim_{r\to 1}(1-r)^{1/q}\log^+M_{\infty}(r,f)=0,$ 
where $M_{\infty}(r,f)=\max_{\vert z\vert\le r}\vert f(z)\vert$ with $0<r<1$.
 Stoll  \cite[Theorem 3.2]{s} proved that the space $F^p$ 
with the topology given by the  family of seminorms 
$\left\{\Vert \cdot\Vert_{q,c}\right\}_{c>0}$ defined for 
$f\in F^q$ as 
  $\Vert f\Vert_{q,c}:=\sum_{n=0}^{\infty}\vert a_n\vert\exp\left(-cn^{1/(q+1)}
\right)<\infty$, becomes a countably normed  Fr\'{e}chet algebra.
It is known that for every $p>1$, $F^p$ is the Fr\'{e}chet envelope of the 
Privalov space $N^p$. 

In this paper, we extend our  study of \cite{m17} 
on the structure of maximal ideals in the algebras $F^p$ ($1<p<\infty$). 
Namely, the obtained characterization of closed maximal 
ideals in $F^p$  from \cite{m17} 
is extended here in terms of topology of uniform convergence
  on compact subsets of $\Bbb D$. 
   \end{abstract}

\maketitle

{\renewcommand{\thefootnote}{}

\vspace{-2mm}

\section{Introduction, Preliminaries and Results}

Let $\Bbb D$ denote the open unit disk 
in the complex plane  and let $\Bbb T$ denote the
boundary of $\Bbb D$.  Let $L^q(\Bbb T)$ $(0<q\le \infty)$ be the
familiar Lebesgue space on the unit circle $\Bbb T$.

The  Privalov class $N^p$ $(1<p<\infty)$ is defined 
as the set of all holomorphic functions $f$ $f$ on $\Bbb D$ such that 
   $$
\sup_{0<r<1}\int_0^{2\pi}(\log^+\vert f(re^{i\theta})\vert)^p\,
\frac{d\theta}{2\pi}<+\infty\eqno(1)
   $$ 
holds, where $\log^+|a|=\max\{\log |a|,0\}$.
 These classes were firstly considered  by Privalov in 
\cite[pages 93--10]{p}, where $N^p$ is denoted as $A_q$. 

Recall  the condition (1) with  $p=1$ defines the 
 Nevanlinna class $N$ of holomorphic functions in $\Bbb D$.
The  Smirnov class $N^+$  is the set
of all functions  $f$ holomorphic on $\Bbb D$ such that
   $$
\lim_{r\rightarrow 1}\int_0^{2\pi}\log^+\vert f(re^{i\theta})\vert
\,\frac{d\theta}{2\pi}=\int_0^{2\pi}\log^+\vert f^*(e^{i\theta})\vert\,
\frac{d\theta}{2\pi}<+\infty,\eqno(2)
   $$
where $f^*$ is the boundary function of $f$ on $\Bbb T$, i.e.,
 $$
 f^*(e^{i\theta})=\lim_{r\rightarrow 1-}f(re^{i\theta})\eqno(3)
  $$ 
is the  radial  limit of a function $f$  which exists for almost 
every $e^{i\theta}\in\Bbb T$. 
 We denote by $H^q$ $(0<q\le\infty)$ the classical  Hardy space on 
$\Bbb D$. 

The following inclusion relations hold true (see \cite{mo, mp1, i1}): 
   $$
N^r\subset N^p\;(r>p),\quad\bigcup_{q>0}H^q\subset
\bigcap_{p>1}N^p,\quad{\rm and}\quad\bigcup_{p>1}N^p\subset N^+
\subset N,\eqno(4)
   $$
where the all  containment relations are proper.
   
The study of the spaces $N^p$ $(1<p<\infty)$ was continued in 1977
by M. Stoll  \cite{s} (with the notation
 $(\log^+H)^\alpha$ in \cite{s}).  Further, the topological and functional 
properties of these spaces  have been  extensively studied by several authors
(see \cite{mo}, \cite{e1},  \cite{e2},  \cite{im} and  
\cite{me4}--\cite{m5}).

M. Stoll \cite[Theorem 4.2]{s}  showed that for every  
 $p>1$ the space $N^p$ (with the noatation $(\log^+H)^{\alpha}$ 
in \cite{s}) equipped with the topology given by the metric 
$d_p$ defined by            
   $$
d_p(f,g)=\Big(\int_0^{2\pi}\big(\log(1+
\vert f^*(e^{i\theta})-g^*(e^{i\theta})\vert)\big)^p\,\frac{d\theta}
{2\pi}\Big)^{1/p},\quad f,g\in N^p,\eqno (5) 
   $$  
becomes an $F$-algebra. This means that $N^p$ is an  $F$-space (a complete 
metrizable topological vector space with the invariant metric)  in which  
multiplication is continuous.

Observe that the function $d_1=d$ defined on the Smirnov class $N^+$
by (5) with $p=1$ induces the metric topology on $N^+$.
N. Yanagihara \cite{y1} proved that 
under this topology, $N^+$ is an $F$-space.

In connection with the spaces $N^p$
$(1<p<\infty)$,  Stoll \cite{s} 
(see also \cite{e1} and \cite[Section 3]{mp3})  also studied the spaces
$F^q$ $(0<q<\infty)$ (with the notation $F_{1/q}$ in \cite{s}), 
consisting of those functions $f$ holomorphic on 
$\Bbb D$ such that    
   $$
 \lim_{r\to 1}(1-r)^{1/q}\log^+M_{\infty}(r,f)=0,\eqno (6) 
   $$ 
where
$$
 M_{\infty}(r,f)=\max_{\vert z\vert\le r}\vert f(z)\vert\quad 
(0<r<1).\eqno (7) 
 $$

In this paper, we will need some Stoll's results
concerning the spaces 
$F^q$ only with $1<q<\infty$. Accordingly, in the sequel, we will assume that 
$q=p>1$ be a fixed real number.

\vspace{2mm}

\noindent{\bf Theorem 1} (see \cite[Theorem 2.2]{s}).   {\it Suppose that
 $f(z)=\sum_{n=0}^\infty a_nz^n$ is a holomorphic function on $\Bbb D$. 
Then the following statements are equivalent:}
\begin{itemize}
\item[(a)] $f\in F^p$;
\item[(b)] {\it there exists a sequence $\{c_n\}_n$ of positive real
numbers with $c_n\to 0$ such that}
   $$
|a_n|\le \exp\left(c_nn^{1/(p+1)}\right),\quad n=0,1,2,\ldots;
\eqno(8) 
   $$
\item[(c)] {\it for any} $c>0$,
          $$
\Vert f\Vert_{p,c}:=\sum_{n=0}^{\infty}\vert a_n\vert\exp\left(-cn^{1/(p+1)}
\right)<\infty.\eqno(9) 
  $$   
\end{itemize}
\vspace{1mm}

\noindent{\it Remark.} Note that in view of the equivalence 
(a)$\Leftrightarrow$(c) of  Theorem 1, 
 by (9) it is well defined  the family of seminorms 
$\left\{\Vert \cdot\Vert_{p,c} \right\}_{c>0}$  on $F^p$.

\vspace{2mm}

Recall that  a locally convex $F$-space is called 
a {\it Fr\'echet space}, and a {\it Fr\'echet algebra} 
is a  Fr\'echet space that is an algebra in which  multiplication
is continuous.  Stoll \cite{s} also proved the following result.

\vspace{2mm}

\noindent{\bf Theorem 2} (see \cite[Theorem 3.2]{s}).   {\it
The space $F^q$ $(0<q<\infty)$ equipped 
with the topology given by the  family of seminorms 
$\left\{\Vert \cdot\Vert_{q,c}\right\}_{c>0}$ defined for 
$f\in F^q$ as 
  $$
\Vert f\Vert_{q,c}:=\sum_{n=0}^{\infty}\vert a_n\vert\exp\left(-cn^{1/(q+1)}
\right)<\infty, \eqno (10) 
  $$
is a countably normed Fr\'{e}chet algebra}.

\vspace{2mm}

Moreover, Stoll \cite{s} defined the family of seminorms  
$\{|\Vert\cdot\Vert_{p,c}|\}_{c>0}$ on $F^p$ given as
   $$
|\Vert f|\Vert_{p,c}=
\int_0^1\exp\left(-c(1-r)^{-1/p}\right)M_p(r,f)\,dr,\quad f\in F^p,\eqno(11)
   $$
where
    $$
M_p(r,f)=\left(\int_0^{2\pi}|f(re^{i\theta})|^p\,\frac{d\theta}{2\pi}\right)
^{1/p}.\eqno(12)
  $$

For our purposes, we will also need the following result.

\vspace{2mm}

\noindent{\bf Theorem 3} (see \cite[Proposition 3.1]{s}).
{\it For each $c>0$, there is a constant $A$ depending only on $p$ and $c$, 
such that
     $$
|\Vert f\Vert|_{p,c}\le \Vert f\Vert_{p,c_1}\quad{\it and}\quad
 \Vert f\Vert_{p,c}\le A |\Vert f\Vert|_{p,c_2},\eqno(13)
      $$
with $c_1=c^{p/(p+1)}$ and $c_2=(c/12)^{p/(p+1)}$. 

Consequently,  $\left\{\Vert \cdot\Vert_{p,c}\right\}_{c>0}$
 and $\left\{|\Vert \cdot\Vert |_{p,c}\right\}_{c>0}$
are equivalent families of seminorms.}
  \vspace{2mm}

It is known  that the Privalov space $N^p$ $(1<p<\infty)$ is not locally 
convex (see \cite[Theorem 4.2]{e1} and \cite[Corollary]{m4}), 
and thus, $N^p$ is properly contained in $F^p$.
Furthermore, $N^p$ is not locally bounded space (see \cite[Theorem 1.1]{mp4}).
Moreover, Stoll proved (\cite[Theorem 4.3]{s}) that for every $p>1$, 
$N^p$ is a dense subspace of $F^p$ and the topology on $F^p$ equiped by the
family of seminorms defined by (10)  is weaker than the topology on  $N^p$ 
induced by the metric $d_p$ defined by (5).
Recall that  Eoff  proved
\cite[Theorem 4.2, the case $p>1$]{e1} that the space $F^p$ is the
{\it Fr\'{e}chet envelope} of $N^p$. For more information on 
the notion of  Fr\'{e}chet envelope, see \cite[Theorem 1]{sh}, 
\cite[Section 1]{m2} and 
  \cite[Corollary 22.3, p. 210]{ke}.

 \vspace{2mm}
\noindent{\it Remark.}  
For $p=1$, the space $F_1$ has been denoted by $F^+$ and has been studied 
by N. Yanagihara in \cite{y2, y1}.  
It was proved in \cite{y2,y1} that $F^+$ is actually 
the containing Fr\'{e}chet space for $N^+$, 
i.e., $N^+$ with the initial topology embeds densely into $F^+$, under
the natural inclusion, and $F^+$ and the Smirnov class $N^+$ have the same 
topological duals.
\vspace{1mm}

Note that the space $F^p$ topologised by the family of
seminorms $\left\{\Vert \cdot\Vert_{p,c}\right\}_{c>0}$ 
given by (10) is metrizable by the metric $\lambda_p$ defined as 
    $ \lambda_p(f,g)=\sum_{n=1}^{\infty}2^{-n}\frac{\Vert f-g\Vert_{p,1/n^
{p/(p+1)}}}{1+\Vert f-g\Vert_{p,1/n^{p/(p +1)}}}$ with $f,g\in F^p.
     $

Since Privalov space $N^p$ and its 
Fr\'{e}chet envelope $F^p$ $(1<p<\infty)$ are algebras, they can be 
also  considered as  rings with respect  to the usual  
ring's operations addition and multiplication. Note that 
these two operations are continuous on the spaces $N^p$ and $F^p$  
in view of the facts that the spaces $N^p$ and $F^p$ are $F$-algebras. 

Motivated by numerous results concerning the ideal structure of some spaces 
of holomorphic functions given in \cite{k1} \cite{mo}, \cite{ma} and 
\cite{be}-\cite{Mor}, 
related investigations on the spaces $N^p$ $(1<p<\infty)$
and their Fr\'{e}chet envelopes were given in 
\cite{mo},  \cite{me4}, \cite{ma}, \cite{mp5}, \cite{mp3}, \cite{m5}
and \cite{m18}.
Notice that a survey of these results was given in \cite{mp6}.
The  $N^p$-analogue of the famous  Beurling's theorem 
for the Hardy spaces $H^q$ $(0<q<\infty)$ \cite{be} was formulated and
 proved in \cite{mp5}.
Moreover, it was showed in \cite[Theorem B]{me4})
that $N^p$ $(1<p<\infty)$ is a ring of Nevanlinna--Smirnov type  
in the sense of Mortini \cite{Mor}. The structure of closed weakly 
dense ideals in $N^p$ was described in \cite{mp3}.
The ideal structure of $N^p$ and the  multiplicative linear functionals  on 
$N^p$ were studied in \cite{mo} and \cite[Theorem ]{m5}. These results are
 similar to those obtained by Roberts and Stoll \cite{RS} for the 
Smirnov class $N^{+}$. 

Motivated by results of Roberts and Stoll given in  \cite[Section 2]{RS2} 
concerning a characterization of  multiplicative linear 
functionals on $F^{+}$ and closed maximal ideals in $F^{+}$,
in \cite{m17} the author of this  paper proved  the analogous 
results for the spaces $F^p$ $(1<p<\infty)$.
These results are  given by  Proposition 5, Proposition 6, 
Theorem 7 and Theorem 8 in \cite{m17}. 
Namely,  if 
 $\lambda\in\Bbb D$  and 
$\gamma_\lambda$ is a functional on $F^p$ defined as
   $$
\gamma_\lambda(f)=f(\lambda)\eqno(14)
   $$ 
for every $f\in F^p$, then by  \cite[Proposition 5]{m17},
 $\gamma_\lambda$ is a continuous multiplicative 
linear functional on $F^p$.
Furthermore, if for any fixed  $\lambda\in\Bbb D$ we define a set 
${\mathcal M}_{\lambda}$  as
  $$
{\mathcal M}_{\lambda}=\{ f\in F^p: f(\lambda)=0\},\eqno(15)
  $$ 
then by \cite[Proposition 6]{m17},
 ${\mathcal M}_{\lambda}$ is a closed maximal ideal in $F^p$. 
Moreover, if $\gamma$ is a nontrivial multiplicative
linear functional on $F^p$,  it is showed  in \cite[Theorem 7]{m17} that  
there exists $\lambda\in \Bbb D$ such that
$$
\gamma(f)=f(\lambda)\eqno(16)
$$ 
for every $f\in F^p$, and in addition,  $\gamma$ is 
a continuous map. Finally, it is proved in \cite[Theorem 8]{m17} 
that every closed maximal ideal in $F^p$ is of the form  
${\mathcal M}_{\lambda}$ for some $\lambda\in \Bbb D$.

Motivated by a result of Igusa \cite{ig},  here we extend
Theorem 8 in   \cite{m17} by the following  result.
  
\vspace{2mm}

\noindent{\bf Theorem 4.} {\it Let $p>1$ and 
let  ${\mathcal M}$ be a maximal  ideal in $F^p$. Then the following 
statements about  ${\mathcal M}$ are equivalent$:$
  \begin{itemize}
\item[(i)] The set ${\mathcal M}$ is closed with respect to the topology of 
uniform convergence on compact subsets of $\Bbb D$$;$ 
   \item[(ii)] $F^p/{\mathcal M}\cong{\bf C}$$;$ 
\item[(iii)]  There exists $\lambda\in \Bbb D$ such that ${\mathcal M}=
{\mathcal M}_
{\lambda}$$;$    
\item[(iv)]  ${\mathcal M}$ is a closed ideal in  $F^p$ with respect 
to the  topology induced on $F^p$ by the family of seminorms 
$\{\Vert\cdot\Vert_{p,c}\}_{c>0}$ defined by $(9)$ and
\item[(v)]  ${\mathcal M}$ is a closed ideal in  $F^p$ with respect 
to the  topology induced on $F^p$ by the family of seminorms 
$\{|\Vert\cdot\Vert|_{p,c}\}_{c>0}$ defined by $(11)$.
 \end{itemize}}

\section{Proof of  Theorem 4}

 In order to obtain a characterization of a maximal ideal 
space of the algebra  $F^p$ with respect to the 
{\it topology of uniform convergence  on compact subsets}  of
$\Bbb D$,  we will need a result of Yanagihara 
 \cite[Lemma 10]{y4} concerning the topological 
 algebras described below.

Let  $A$ be a topological algebra  over the field $\Bbb C$ with identity
1, locally convex and  commutative. The topology  of $A$ 
is defined by a countable  family  of seminorms 
$\{\Vert\cdot\Vert_{\alpha}\}_{\alpha\in S}$ 
($S=\Bbb N$ or $S$ is a finite subset of the set
of positive integers  $\Bbb  N$) 
for which  $\Vert 1\Vert_{\alpha}=1$ and for $a,b\in A$ there holds
      $$
\Vert ab\Vert_{\alpha}\le\Vert a\Vert_{\alpha}\cdot\Vert b\Vert_{\alpha}
\quad {\rm for\;\;all}\;\;a,b\in A,\;\alpha\in S.\eqno(17)
      $$
For an  $\alpha\in S$ let $E_{\alpha}=\{a\in A:\,\Vert a\Vert_{\alpha}=0\}$.
$E_{\alpha}$ is obviously an ideal of the algebra $A$. 
For any $a\in A$ we define the coset  $\tilde{a}$ 
as $\tilde{a}=a+E_{\alpha}\in A/{E_{\alpha}}$.
The  quotient space $A/{E_{\alpha}}$ is a normed space
with the assciated norm 
$\Vert\tilde{a}\Vert_{\alpha}= \Vert a\Vert_{\alpha}$, 
$\alpha\in S$. By (17) we have 
       $$
\Vert\tilde{a}\tilde{b}\Vert_{\alpha}\le\Vert\tilde{a}\Vert_{\alpha}
\cdot \Vert \tilde{b}\Vert_{\alpha}
\quad {\rm for\,\, all}\;\;\tilde{a},\tilde{b}\in A/{E_{\alpha}}
  ,\,\,\alpha\in S.\eqno(18)
      $$
The completion of  the space  $A/{E_{\alpha}}$ with respect to the 
norm  $\Vert\cdot \Vert_{\alpha}$ is denoted by $A_{\alpha}^*$. 
Then we have the following result.

\vspace{2mm}

\noindent{\bf Lemma 5.} (\cite[Lemma 10]{y4}). 
{\it Let  $A$ be a topological algebra  described above.
Then for any  $f\in A$ there exists a complex number ${\lambda}_f$  
 depending on $f$ such that ${\lambda}_f-f=({\lambda}_f1-f)$ is not
 invertible element of the algebra $A$.}
 
\vspace{2mm}

\noindent{\it Proof of Theorem 4.}  
(i)$\Rightarrow$(ii). 
Obviously, there holds  $F^p/{{\mathcal M}}\supseteq \Bbb C$. 
For  $f\in F^p$ denote
    $$
[f]=f+{\mathcal M}\in F^p/{{\mathcal M}}.\eqno(19)
    $$
Now we define the family of seminorms $\{\Vert\cdot\Vert_r\}$ $(0\le r<1)$ in  
$F^p/{{\mathcal M}}$ as follows: 
     $$
\Vert [f]\Vert_r=\inf_{h\in {\mathcal M}}(\max_{|z|=r}|f(z)+h(z)|),\quad 0
\le r<1.\eqno(20)
     $$
Then obviously, we have
     $$
\Vert [fg]\Vert_r\le\Vert [f]\Vert_r\cdot \Vert [g]\Vert_r,\quad 0\le r<1.\eqno(21)
     $$
By Lemma 5, to each $[f]\in F^p/{{\mathcal M}}$ corresponds 
a number $\lambda\in \Bbb C$ such that $\lambda-[f]$ 
is not invertible element of  $F^p/{{\mathcal M}}$. 
However, by the maximality of the ideal  ${\mathcal M}$,
we conclude that  $F^p/{{\mathcal M}}$ is a field, 
and thus, $\lambda-f$ must belong to  ${\mathcal M}$, 
i.e., $\lambda\in[f]$. Hence, we obtain 
(cf. \cite[Proof of Theorem 8]{m17} with application of Arens's 
result \cite{ar})
     $$
F^p/{{\mathcal M}}\cong\Bbb C.\eqno(22)
     $$
\indent (ii)$\Rightarrow$(iii). Let  $\lambda$ be in  the coset $[z]\in 
F^p/{{\mathcal M}}$. Then $z- \lambda\in {\mathcal M}$, and hence,  
$\lambda\in \Bbb D$.
\\\indent
For each $f\in F^p$ we have
     $$
f(z)-f(\lambda)=A(z)(z-\lambda).\eqno(23)
     $$
It is easy to see that $A(z)\in F^p$, and thus, 
$f(z)-f(\lambda)\in {\mathcal M}$.
Therefore, if  $f(z)\in {\mathcal M}$, then $f(\lambda)\in {\mathcal M}$, 
whence it follows that $f(\lambda)=0$. This shows that 
${\mathcal M}={\mathcal M}_{\lambda}$.\\\indent
(iii)$\Rightarrow$(i). This implication is evident from the  theorem of
 Hurwitz.\\
\indent (iii)$\Rightarrow$(iv). This implication immediately follows 
from  \cite[Proposition 6]{m17}.\\\indent
(iv)$\Rightarrow$(iii). This implication immediately follows 
from  \cite[Theorem 8]{m17}.\\\indent
(iv)$\Leftrightarrow$(v). This equivalence immediately follows 
from  Theorem 3.
This completes the proof of Theorem 4.

\end{document}